\documentclass[11pt]{article}
\usepackage{cite}
\usepackage{mathrsfs}
\usepackage{amssymb,amsfonts,amsmath}
\usepackage{enumerate}
\usepackage{indentfirst}
\usepackage{latexsym,bm}
\usepackage[noend]{algpseudocode}
\usepackage{algorithmicx,algorithm}
\usepackage{algorithm}
\usepackage{makeidx}
\usepackage{fancybox}
\usepackage{color}
\usepackage{multicol}
\usepackage[latin1]{inputenc}
\usepackage{graphicx}
\usepackage{epstopdf}
\usepackage{pstricks,pst-node,pst-text,pst-3d}
\usepackage{tikz}
\newtheorem{thm}{Theorem}[section]

\newtheorem{pro}[thm]{Proposition}

\newenvironment{pf}{{\noindent \it \bf Proof:}}{{\hfill$\Box$}\\}

\begin{document}

\title{\bf Strong subgraph 2-arc-connectivity and arc-strong connectivity of Cartesian product of digraphs}
\author{Yiling Dong{ }$^{1}$, Gregory Gutin$^{2}$ and Yuefang Sun$^{3,}$\footnote{Corresponding author. Yuefang Sun was supported by Yongjiang Talent Introduction Programme of Ningbo and Zhejiang Provincial Natural Science Foundation of China under Grant No. LY20A010013.} \\
$^{1}$ School of Mathematics and Statistics,
Ningbo University\\
Zhejiang 315211, P. R. China, dongyilingnpc@163.com\\
$^{2}$ Department of Computer Science\\ Royal Holloway University of London\\
Egham, Surrey, UK, gutin@cs.rhul.ac.uk\\
$^{3}$ School of Mathematics and Statistics,
Ningbo University\\
Zhejiang 315211, P. R. China, sunyuefang@nbu.edu.cn}
\date{}
\maketitle

\begin{abstract}
Let $D=(V,A)$ be a digraph of order $n$, $S$ a subset of $V$ of size $k$ and $2\le k\leq n$. A strong subgraph $H$ of $D$ is called an $S$-strong subgraph if $S\subseteq V(H)$. A pair of $S$-strong subgraphs $D_1$ and $D_2$ are said to be  arc-disjoint if $A(D_1)\cap A(D_2)=\emptyset$.
Let $\lambda_S(D)$ be the maximum number of arc-disjoint $S$-strong subgraphs in $D$. The strong subgraph $k$-arc-connectivity is defined as $$\lambda_k(D)=\min\{\lambda_S(D)\mid S\subseteq V(D), |S|=k\}.$$ The parameter $\lambda_k(D)$ can be seen as a
generalization of classical edge-connectivity of undirected graphs.

In this paper, we first obtain a formula for the arc-connectivity of Cartesian product $\lambda(G\Box H)$ of two digraphs $G$ and $H$ generalizing a formula for edge-connectivity of  Cartesian product of two undirected graphs obtained by Xu and Yang (2006). Then we study the strong subgraph 2-arc-connectivity of Cartesian product $\lambda_2(G\Box H)$ and prove that $ \min\left \{ \lambda \left ( G \right ) \left | H \right | , \lambda \left ( H \right ) \left |G \right |,\delta ^{+  } \left ( G \right )+ \delta ^{+  } \left ( H \right ),\delta ^{-   } \left ( G \right )+ \delta ^{-   } \left ( H \right )  \right \}\ge\lambda_2(G\Box H)\ge \lambda_2(G)+\lambda_2(H)-1.$ The upper bound for  $\lambda_2(G\Box H)$ is sharp and is a simple corollary of the formula for $\lambda(G\Box H)$.  
The lower bound for $\lambda_2(G\Box H)$ is either sharp or almost sharp i.e. differs by 1 from the sharp bound. We also obtain exact values for $\lambda_2(G\Box H)$, where $G$ and $H$ are digraphs from some digraph families. 
\vspace{0.3cm}
\vspace{0.3cm}
\end{abstract}


\section{Introduction}\label{sec:intro}

We refer the readers to \cite{Bang-Jensen-Gutin, Bondy} for graph theoretical notation and terminology not given here. Note that all digraphs considered in this paper have no parallel arcs or loops. A digraph $D$ is {\em symmetric} if it can be obtained from
its underlying undirected graph $G$ by replacing each edge of $G$
with the corresponding arcs of both directions, that is,
$D=\overleftrightarrow{G}$. The order $|G|$ of a (di)graph $G$ is the number of vertices in $G.$ 
Let $\overleftrightarrow{T}_n$ be the symmetric digraph whose underlying undirected graph is a tree of order $n$. We use $\overrightarrow{C}_n$ and $\overleftrightarrow{K}_n$ to denote the cycle and complete digraph of order $n$, respectively.

For a graph $G=(V,E)$ and a set $S\subseteq V$ of at least two vertices,
an {\em $S$-Steiner tree} or, simply, an {\em $S$-tree} is a subgraph
$T$ of $G$ which is a tree with $S\subseteq V(T)$. Two
$S$-trees $T_1$ and $T_2$ are said to be {\em edge-disjoint} if
$E(T_1)\cap E(T_2)=\emptyset$. Two arc-disjoint $S$-trees $T_1$ and $T_2$ are said to be {\em internally disjoint} if $V(T_1)\cap V(T_2)=S$. The {\em
generalized local connectivity} $\kappa_S(G)$ is the maximum number
of internally disjoint $S$-trees in $G$. For an integer $k$ with
$2\leq k\leq n$, the {\em generalized $k$-connectivity} \cite{Hager} is defined as
$$\kappa_k(G)=\min\{\kappa_S(G)\mid S\subseteq V(G), |S|=k\}.$$
Similarly, the {\em generalized local edge-connectivity} $\lambda_S(G)$ is the maximum number of edge-disjoint $S$-trees in $G$. For an integer $k$ with $2\leq k\leq n$, the {\em generalized $k$-edge-connectivity} \cite{Li-Mao-Sun} is defined as
$$\lambda_k(G)=\min\{\lambda_S(G)\mid S\subseteq V(G), |S|=k\}.$$
Let $\kappa(G)$ and $\lambda(G)$ denote the classical vertex-connectivity and edge-connectivity of an undirected graph $G.$
Observe that $\kappa_2(G)=\kappa(G)$ and $\lambda_2(G)=\lambda(G)$, hence, these two parameters are generalizations of classical connectivity of undirected graphs and are also called tree connectivity. Now the topic of tree connectivity has become an established area in graph theory, see a recent monograph \cite{Li-Mao5} by Li and Mao on this topic.

To extend generalized $k$-connectivity to directed graphs, Sun,
Gutin, Yeo and Zhang \cite{Sun-Gutin-Yeo-Zhang} observed that in the
definition of $\kappa_S(G)$, one can replace ``an $S$-tree'' by ``a
connected subgraph of $G$ containing $S$.'' Therefore, they defined {\em strong subgraph $k$-connectivity} by replacing ``connected'' with ``strongly connected'' (or, simply, ``strong'') as follows. Let $D=(V,A)$ be a digraph of order $n$, $S$ a subset of $V$ of size $k$ and $2\le
k\leq n$. An {\em S-strong subgraph} is a strong subgraph $H$ of $D$ such that $S\subseteq V(H)$. $S$-strong subgraphs $D_1, \dots , D_p$ are
said to be {\em internally disjoint} if $V(D_i)\cap V(D_j)=S$ and
$A(D_i)\cap A(D_j)=\emptyset$ for all $1\le i<j\le p$. Let
$\kappa_S(D)$ be the maximum number of internally disjoint $S$-strong
digraphs in $D$. The {\em strong subgraph
$k$-connectivity} \cite{Sun-Gutin-Yeo-Zhang} is defined as
$$\kappa_k(D)=\min\{\kappa_S(D)\mid S\subseteq V, |S|=k\}.$$

As a natural counterpart of the strong subgraph $k$-connectivity, Sun and Gutin \cite{Sun-Gutin} introduced the concept of strong subgraph $k$-arc-connectivity.
Let $D=(V(D),A(D))$ be a digraph of order $n$, $S\subseteq V$ a
$k$-subset of $V(D)$ and $2\le k\leq n$.
Let $\lambda_S(D)$ be the maximum number of arc-disjoint $S$-strong digraphs in $D$. The {\em strong subgraph
$k$-arc-connectivity} is defined as
$$\lambda_k(D)=\min\{\lambda_S(D)\mid S\subseteq V(D), |S|=k\}.$$ 
Note that $\kappa_k(D)$ and $\lambda_k(D)$ are not only natural
extensions of tree connectivity, but also could be seen as
generalizations of connectivity and edge-connectivity of undirected
graphs as $\kappa_2(\overleftrightarrow{G})=\kappa(G)$
\cite{Sun-Gutin-Yeo-Zhang} and
$\lambda_2(\overleftrightarrow{G})=\lambda(G)$ \cite{Sun-Gutin}.
For more information on the topic of strong subgraph connectivity of digraphs, the readers can see \cite{Sun-Gutin2} for a recent survey.

In this paper, we continue research on strong subgraph  arc-connectivity and focus on the strong subgraph 2-arc-connectivity of Cartesian products of digraphs. It is well known that Cartesian products of digraphs are of interest in graph theory and its applications; see a recent survey chapter by Hammack \cite{Hammack} considering many results on Cartesian products of digraphs.

In the next section we introduce terminology and notation on Cartesian products of digraphs and give a simple yet useful upper bound on $\lambda_2(D),$ where $D$ is Cartesian product of any digraphs $G$ and $H$ i.e. $D=G\Box H$.

In Section \ref{sec:exect}, we prove that 
$$\lambda  \left ( G \Box H\right )= \min\left \{ \lambda \left ( G \right ) \left | H \right | , \lambda \left ( H \right ) \left |G \right |,\delta ^{+  } \left ( G \right )+ \delta ^{+  } \left ( H \right ),\delta ^{-   } \left ( G \right )+ \delta ^{-   } \left ( H \right )  \right \}$$ for every pair  $G$ and $H$ of strong digraphs, each of order at least 2.\footnote{Note that the case of at least one of two digraphs having just one vertex in $\lambda  \left ( G \Box H\right )$ is trivial. Thus, we will henceforth assume  that each of the two digraphs is of order at least 2. The same holds for $\lambda_2  \left ( G \Box H\right )$. }

In Section~\ref{sec:1product} we prove that $$ \min\left \{ \lambda \left ( G \right ) \left | H \right | , \lambda \left ( H \right ) \left |G \right |,\delta ^{+  } \left ( G \right )+ \delta ^{+  } \left ( H \right ),\delta ^{-   } \left ( G \right )+ \delta ^{-   } \left ( H \right )  \right \}$$ and $\lambda_2(G)+\lambda_2(H)-1$ are an upper bound and a lower bound, respectively, for 
$\lambda_2(G\Box H)$. The upper bound follows from the formula for $\lambda  \left ( G \Box H\right )$ and thus it is tight.
Unfortunately, we do not know whether this lower bound is tight or not, but by Theorem \ref{thmd1} (mentioned below), the gap with a tight bound is at most 1. 

In Section~\ref{sec:product}, we obtain exact values for the strong subgraph 2-arc-connectivity of Cartesian products of some digraph classes; our results are collated in Theorem~\ref{thmd1}. For the classes of strong digraphs considered in Theorem \ref{thmd1}, we have $\lambda_2(G\Box H)= \lambda_2(G) +\lambda_2(H).$ 



\section{Cartesian product of digraphs}\label{sec:cp}

For a positive integer $n$, let $[n]=\{1,2,\dots ,n\}.$

Let $G$ and $H$ be two digraphs with $V(G)=\{u_i \mid 1\leq i\leq
n\}$ and $V(H)=\{v_j \mid 1\leq j\leq m\}$. The {\em Cartesian product} $G\Box H$ of two digraphs $G$ and $H$ is a
digraph with vertex set
$$V(G\Box H)=V(G)\times V(H)=\{(x, x')\mid x\in V(G), x'\in V(H)\}$$
and arc set $$A(G\Box H)=\{(x,x')(y,y')\mid xy\in A(G),~x'=y',~or~x=y,~x'y'\in A(H)\}.$$~We will use $u_{i,j}$ to denote $(u_i,v_j)$ in the rest of the paper.~By definition, we know the
Cartesian product is associative and commutative, and $G\Box H$ is
strongly connected if and only if both $G$ and $H$ are strongly
connected \cite{Hammack}.

\begin{figure}[htbp]
\small
\centering
  \includegraphics[width=9cm]{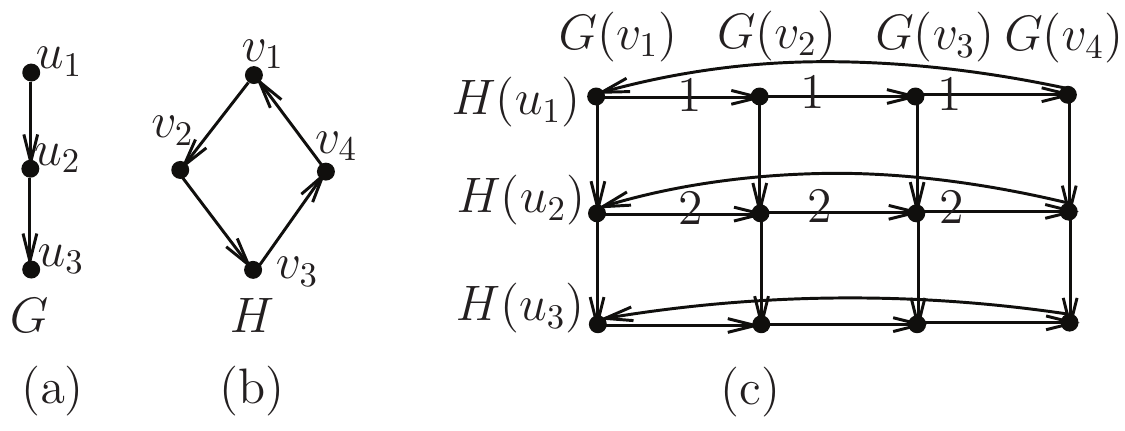}
\caption{Two digraphs $G$, $H$ and their Cartesian
product.}
\label{figure1}
\end{figure}

We use $G(v_j)$ to
denote the subgraph of $G\Box H$ induced by vertex set
$\{u_{i,j}\mid 1\leq i\leq n\}$ where $1\leq j\leq m$,~and use
$H(u_i)$ to denote the subgraph of $G\Box H$ induced by vertex set
$\{u_{i,j}\mid 1\leq j\leq m\}$ where $1\leq i\leq n$.~Clearly,~we
have $G(v_j)\cong G$ and $H(u_i)\cong H$. (For example,~as shown in
Fig. \ref{figure1},~$G(v_j)\cong G$ for $1\leq j\leq 4$ and
$H(u_i)\cong H$ for $1\leq i\leq 3$).~For $1\leq j_1\neq j_2\leq m$,~
the vertices $u_{i,j_1}$ and $u_{i,j_2}$ belong to the same
digraph $H(u_i)$ where $u_i\in V(G)$;~we call $u_{i,j_2}$ the
{\em vertex corresponding to} $u_{i,j_1}$ in $G(v_{j_2})$;~for
$1\leq i_1\neq i_2\leq n$, we call $u_{i_2,j}$ the vertex
corresponding to $u_{i_1,j}$ in $H(u_{i_2})$.~Similarly,~we can
define the subgraph {\em corresponding} to some subgraph.~For example,~
in the digraph (c) of Fig. \ref{figure1},~let $P_1$~$(P_2)$ be the
path labelled 1 (2) in $H(u_1)~(H(u_2))$, then $P_2$ is called the
path {\em corresponding} to $P_1$ in $H(u_2)$.

It follows from the definition of strong subgraph 2-arc-connectivity  that for any digraph $D$, $\lambda_2(D)\le \min\{\delta^+(D), \delta^-(D)\}$ \cite{Sun-Gutin}. We will use this inequality 
in Section~\ref{sec:product}. Note that if $D=G\Box H$ then $\delta^+(D)=\delta^+(G)+\delta^+(H)$
and $\delta^-(D)=\delta^-(G)+\delta^-(H).$



\section{Formula for arc-connectivity of Cartesian product of two digraphs}\label{sec:exect}

Xu and Yang \cite{Xu-Yang} (see also \cite{Spa} and \cite[Theorem 5.5]{Imrich-Klavzar-Rall}) proved that 
\begin{equation}\label{eq3}
\lambda(G\Box H)=\min\{\lambda(G)|V(H)|,\lambda(H)|V(G)|,\delta(G)+\delta(H)\}  
\end{equation}
for all connected undirected graphs $G$ and $H,$ each with at least two vertices.
Since $\lambda(\overleftrightarrow{Q})=\lambda(Q)$ for every undirected graph $Q$, Formula (\ref{eq3}) can be easily extended 
to symmetric digraphs. 
In this section, we generalise Formula (\ref{eq3}) to all strong digraphs. 

Clearly, $\lambda\left (  D  \right ) \le \min\left \{\delta ^{+  } \left ( D \right ),\delta ^{-   } \left ( D \right )  \right \} $ for every digraph $D$. Hence, for any two strong digraphs $G$ and $H$, we have
\begin{equation}\label{2}
\begin{split}
\lambda\left (  G \Box H \right ) \le \min\left \{\delta ^{+  } \left ( G \Box H \right ),\delta ^{-   } \left ( G \Box H \right )  \right \}  
\\=\min\left \{ \delta ^{+  } \left ( G \right )+ \delta ^{+  } \left ( H \right ),\delta ^{-   } \left ( G \right )+ \delta ^{-   } \left ( H \right )   \right \}.
\end{split}
\end{equation}

Furthermore, by the definitions of arc-strong connectivity and Cartesian product of digraphs, we have
\begin{equation}\label{3}
\lambda  \left ( G \Box H\right )\le  \lambda\left ( G \right ) \left | H \right | 
\end{equation}
and
\begin{equation}\label{4}
\lambda  \left ( G \Box H\right )\le  \lambda\left ( H \right ) \left | G \right |.
\end{equation}

The inequalities (\ref{2}), (\ref{3}) and (\ref{4}) imply that $$\lambda  \left ( G \Box H\right )\le \min\left \{ \lambda \left ( G \right ) \left | H \right |, \lambda \left ( H \right ) \left |G \right |,\delta ^{+  } \left ( G \right )+ \delta ^{+  } \left ( H \right ),\delta ^{-   } \left ( G \right )+ \delta ^{-   } \left ( H \right )  \right \}.$$ In fact, we can furthermore prove that the equality holds and it could be seen as a digraph extension of (\ref{eq3}).

\begin{thm}\label{arc-connectivity}
Let $G$ and $H$ be two strong digraphs, each of order at least 2. Then
$$\lambda  \left ( G \Box H\right )= \min\left \{ \lambda \left ( G \right ) \left | H \right | , \lambda \left ( H \right ) \left |G \right |,\delta ^{+  } \left ( G \right )+ \delta ^{+  } \left ( H \right ),\delta ^{-   } \left ( G \right )+ \delta ^{-   } \left ( H \right )  \right \}.$$
\end{thm}
\begin{pf}
Let ${S}\subseteq A\left (G \Box H  \right ) $ be an arc-cut set of  $G \Box H$  with $\left|{S}  \right | =\lambda  \left ( G \Box H\right )$. It suffices to show that 
$$|S|\ge \min\left \{ \lambda \left ( G \right ) \left | H \right | , \lambda \left ( H \right ) \left |G \right |,\delta ^{+  } \left ( G \right )+ \delta ^{+  } \left ( H \right ),\delta ^{-   } \left ( G \right )+ \delta ^{-   } \left ( H \right )  \right \}.$$

If $\left |{S}  \right |\ge  \min\left \{ \lambda \left ( G \right ) \left | H \right | , \lambda \left ( H \right ) \left |G \right | \right \}$, then the inequality clearly holds. 

Therefore, we assume that $\left |{S}  \right |<  \min\left \{ \lambda \left ( G \right ) \left | H \right | , \lambda \left ( H \right ) \left |G \right | \right \}$ in the following argument and in this case it suffices to show that 
$\left | {S} \right | \ge\delta ^{+  } \left ( G \Box H \right )$ or $\left | {S} \right | \ge\delta ^{-} \left ( G \Box H \right )$.

Now there must exist a strong component $B$ of $G \Box H-{S}$ which contains some $G\left ( v_{j}\right )$, say $G\left ( v_{1}\right )$, (as $\left | {S} \right | <\lambda \left ( G \right ) \left | H \right |$) and some $H\left ( u_{i} \right )$, say $H\left ( u_{1} \right )$,  in
$G \Box H-{S}$ (as $\left | {S} \right | <\lambda \left ( H \right ) \left | G \right |$). Let $\left ( u,v \right )\in V(G \Box H)\setminus V(B)$.

We want to prove that $\left | {S} \right | \ge d^{+  } \left (( u,v )\right )$ by the following operation that assigns each out-neighbor of $\left ( u,v \right )$ in $G\Box H$ a unique arc from ${S}$:

We first consider out-neighbors of $\left ( u, v \right )$ in $G\left ( v \right )$. Let $\left ( {u}', v \right )$ be an out-neighbor of $\left ( u, v \right )$ in $G\left ( v \right )$. If the arc $a=\left ( u ,v  \right )\left ( {u}'  ,v  \right )\in {S}$, we assign $a$ to $\left ( {u}', v \right )$. Otherwise, we must have $\left ( {u}', v  \right )\notin B$. Therefore, the subdigraph of $G \Box H-{S}'$ induced by $V(H\left ( {u}' \right ))$ is not strong and so $H\left ( {u}' \right )$  contains at least one arc from ${S}$, and we assign this arc to $\left ( {u}'  ,v  \right )$.

We next consider out-neighbors of $\left ( u, v \right )$ in $H\left ( u \right )$. Let $\left ( u, {v}' \right )$ be an out-neighbor of $\left ( u, v \right )$ in $H\left ( u \right )$. If ${a}'=\left ( u, v \right )\left ( u  ,{v}'  \right )\in {S}$, we assign ${a}'$ to $\left ( u, {v}' \right )$. Otherwise, we must have $\left ( u, {v}'  \right )\notin B$. Therefore, the subdigraph of $G \Box H-{S}'$ induced by $V(G\left ( {v}' \right ))$ is not strong and so $G\left ( {v}' \right )$  contains at least one arc from ${S}$, and we assign this arc to $\left ( u, {v}'  \right )$.

The above operations mean that $\left | {S} \right | \ge d^{+  } \left (( u,v )\right )  \ge \delta ^{+  } \left ( G \Box H \right )$. With a similar argument, we can prove that $\left | {S} \right | \ge\delta ^{-  } \left ( G \Box H \right )$. This completes the proof.
\end{pf}

\section{General bounds}\label{sec:1product}

By Theorems~\ref{thmd1} and ~\ref{arc-connectivity}, and the fact that   
$\lambda _{k} \left ( D \right )\le \lambda\left (  D  \right )$  for any digraph $D$\cite{Sun-Gutin}, we have the following sharp upper bound for $\lambda _{2} \left ( G \Box H\right )$. 

\begin{thm}
Let $G$ and $H$ be two strong digraphs, each with at least two vertices. Then
$$\lambda _{2} \left ( G \Box H\right )\le \min\left \{ \lambda \left ( G \right ) \left | H \right | , \lambda \left ( H \right ) \left |G \right |,\delta ^{+  } \left ( G \right )+ \delta ^{+  } \left ( H \right ),\delta ^{-   } \left ( G \right )+ \delta ^{-   } \left ( H \right )  \right \}.$$ Moreover, this bound is sharp.
\end{thm}

Now we will provide a lower bound for $\lambda_2(G\Box H)$ for strong digraphs $G$ and $H$.

\begin{thm}\label{thmd}
Let $G$ and $H$ be two strong digraphs. We have $$\lambda_2(G\Box H)\geq
\lambda_2(G)+ \lambda_2(H)-1.$$ 
\end{thm}

\begin{pf}
It suffices to show that there are at least $\lambda_2(G)+
\lambda_2(H)-1 $ arc-disjoint $S$-strong subgraphs for any
$S\subseteq V(G\Box H)$ with $|S|=2$. Let $S=\{x, y\}$ and 
consider the following two cases.

{\em Case 1}: $x$ and $y$ are in the same $H(u_i)$ or $G(v_j)$ for
some $1\leq i\leq n, 1\leq j\leq m$. We will prove that, in this case, $\lambda_2(G\Box H)\geq \lambda_2(G)+ \lambda_2(H).$
Without loss of generality, we may
assume that $x=u_{1,1},~y=u_{1,2}$. We know there are at least
$\lambda_2(H)$ arc-disjoint $S$-strong subgraphs in the
subgraph $H(u_1)$, and so it suffices to find the remaining $\lambda_2(G)$ $S$-strong subgraphs in $G\Box H$.

We know there are at least $\lambda_2(G)$ arc-disjoint $\{x, u_{2,1}\}$-strong subgraphs, say $D_i(v_1)~(i\in [\lambda_2(G)])$, in $G(v_1)$. For each $i\in
[\lambda_2(G)]$, we can choose an out-neighbor, say
$u_{t_i,1}$~$(i\in [\lambda_2(G)])$, of $x$ in $D_i(v_1)$ such
that these out-neighbors are distinct. Then in $H(u_{t_i})$, we know
there are $\lambda_2(H)$ arc-disjoint $\{u_{t_i,1}, u_{t_i,2}\}$-strong subgraphs, we choose one such strong subgraph, say $D(H(u_{t_i}))$.~For each $i\in [\lambda_2(G)]$,~let $D_i(v_2)$ be the $\{u_{t_i,2}, y\}$-strong subgraph corresponding to $D_i(v_1)$
in $G(v_2)$. We now construct the remaining $\lambda_2(G)$ $S$-strong
subgraphs by letting $D_i=D_i(v_1)\cup
D(H(u_{t_i}))\cup D_i(v_2)$ for each $i\in [\lambda_2(G)]$.
Combining the former $\lambda_2(H)$ arc-disjoint $S$-strong subgraphs with the $\lambda_2(G)$ $S$-strong
subgraphs, we can obtain $\lambda_2(G)+ \lambda_2(H)$ strong subgraphs. Observe all these strong subgraphs are arc-disjoint.


{\em Case 2.} $x$ and $y$ belong to distinct $H(u_i)$ and $G(v_j)$.
Without loss of generality, we may assume that $x=u_{1,1},~y=u_{2,2}$.

There are at least $\lambda_2(G)$ arc-disjoint $\{x,
u_{2,1}\}$-strong subgraphs, say $D_i(v_1) $ $(i\in [\lambda_2(G)])$, in $G(v_1)$.~For each $i\in [\lambda_2(G)]$,~we can
choose an out-neighbor,~say $u_{t_i,1}$~$(i\in [\lambda_2(G)])$,
of $x$ in $D_i(v_1)$ such that these out-neighbors are distinct.
Then in $H(u_{t_i})$, we know that there are $\lambda_2(H)$ arc-disjoint
$\{u_{t_i,1}, u_{t_i,2}\}$-strong subgraphs; we choose one such strong subgraph,~say $D(H(u_{t_i}))$. For each $i\in [\lambda_2(G)]$,~let $D_i(v_2)$ be the $\{u_{t_i,2},~y\}$-strong subgraph corresponding to $D_i(v_1)$ in $G(v_2)$. We now construct the
$\lambda_2(G)$ $S$-strong subgraphs by letting
$D_i=D_i(v_1)\cup D(H(u_{t_i}))\cup D_i(v_2)$ for each $i\in
[\lambda_2(G)]$.

Similarly, there are at least $\lambda_2(H)$ arc-disjoint $\{x, u_{1,2}\}$-strong subgraphs, say $D'_j(u_1)~(j\in [\lambda_2(H)])$, in $H(u_1)$. For each $j\in[\lambda_2(H)]$, 
we can choose an out-neighbor, 
say $u_{1,t'_j}$~$(j\in [\lambda_2(H)])$, of $x$ in $D'_j(u_1)$ such that these out-neighbors are distinct. Then in $G(v_{t'_j})$, we
know there are $\lambda_2(G)$ arc-disjoint $\{u_{1,t'_j}, u_{2,t'_j}\}$-strong subgraphs, we choose one such strong subgraph, say $D(G(v_{t'_j}))$. For each $j\in [\lambda_2(H)]$, let $D'_j(u_2)$ be the $\{u_{2,t'_j}, y\}$-strong subgraph corresponding
to $D'_j(u_1)$ in $H(u_2)$. We now construct the other
$\lambda_2(H)$ $S$-strong subgraphs by letting
$D'_j=D'_j(u_1)\cup D(G(v_{t'_j}))\cup D'_j(u_2)$ for each $j\in
[\lambda_2(H)]$.

{\em Subcase 2.1.} $t_i\neq 2$ for any $i\in [\lambda_2(G)]$ and
$t'_j\neq 2$ for any $j\in [\lambda_2(H)]$, that is,~$u_{2,1}$ was
not chosen as an out-neighbor of $u_{1,1}$ in $G(v_1)$ and
$u_{1,2}$ was not chosen as an out-neighbor of $u_{1,1}$ in
$H(u_1)$.~We can check the above $\lambda_2(G)+ \lambda_2(H)$ strong
subgraphs are arc-disjoint.

{\em Subcase 2.2.} $t_i=2$ for some $i\in [\lambda_2(G)]$ or
$t'_j=2$ for some $j\in [\lambda_2(H)]$,~that is,~$u_{2,1}$ was
chosen as an out-neighbor of $u_{1,1}$ in $G(v_1)$ or $u_{1,2}$
was chosen as an out-neighbor of $u_{1,1}$ in $H(u_1)$. 
Without loss of generality,~we may assume that $t_i=2$ and $t'_j \neq 2 $,~that is,~$u_{2,1}$ was chosen as an out-neighbor of $u_{1,1}$ in $G(v_1)$ and $u_{1,2}$ was not chosen as an out-neighbor of $u_{1,1}$ in
$H(u_1)$.~When $A(D_i)\cap A(D_j) \neq \emptyset$,~we can get $\lambda_2(G)+ \lambda_2(H)-1$ arc-disjoint $S$-strong subgraphs.~Otherwise,~we can check the above $\lambda_2(G)+ \lambda_2(H)$ strong subgraphs are arc-disjoint and get the desired $S$-strong subgraphs.

{\em Subcase 2.3.}  $t_i=2$ for some $i\in [\lambda_2(G)]$ and
$t'_j=2$ for some $j\in [\lambda_2(H)]$,  we replace
$D_1$,~$D'_1$ by $\overline{D}_1$, $\overline{D'}_1$, respectively
as follows: let $\overline{D}_1= D_1(v_1)\cup D(H(u_{t_1}))$ and
$\overline{D'}_1= D'_1(u_1)\cup D_1(v_2)$. We can check that the
current $\lambda_2(G)+ \lambda_2(H)$ strong subgraphs are
arc-disjoint.

Hence, the bound holds. This completes the proof.
\end{pf}

\section{Exact values for digraph classes}\label{sec:product}

In this section, we will obtain exact
values for the strong subgraph 2-arc-connectivity of Cartesian
product of two digraphs belonging to some digraph classes.

\begin{pro}\label{p1}
We have $ \lambda_2(\overrightarrow{C}_n\Box \overrightarrow{C}_m)=2. $
\end{pro}
\begin{figure}[htbp]
\small
\centering
  \includegraphics[width=9cm]{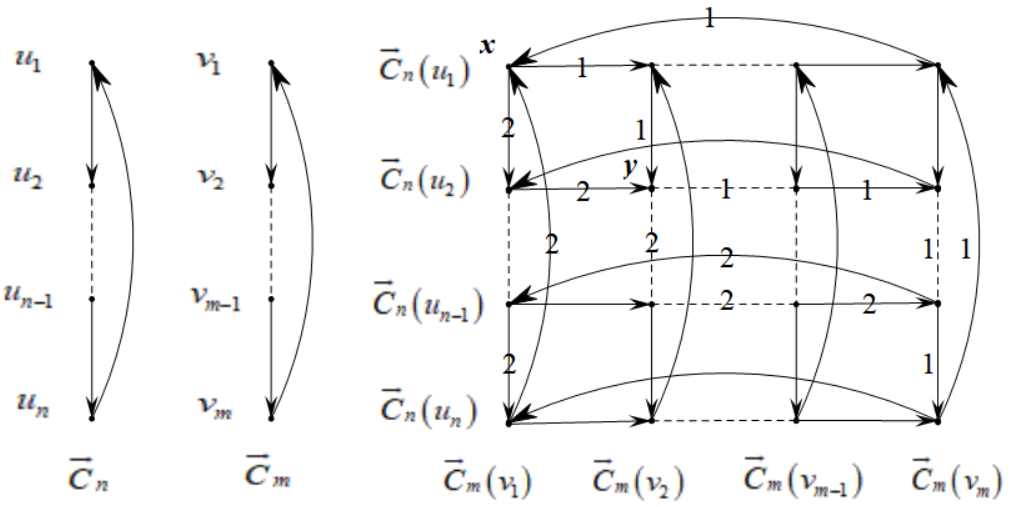}
\caption{Cartesian product of two dicycles.}
\label{figure3.1}
\end{figure}

\begin{pf}\,   
Let $S=\left \{ x, y \right \} $,~we just consider the case that $ x,~y$  are  neither in the same $ \overrightarrow{C}_n\left ( u_{i}  \right )  $ nor in the same  $ \overrightarrow{C}_m\left ( v_{j}  \right )  $ for some $ 1\leq i \leq n $,~$ 1\leq j \leq m $, since the arguments for remaining cases are similar.~Without loss of generality, we may assume that  $ x= u_{1,1},~y=u_{2,2} $.~We can  get  two arc-disjoint  $S$-strong subgraphs  in $\overrightarrow{C}_n\Box \overrightarrow{C}_m$,~say $ D_{1} $ and $ D_{2} $ (as shown in Fig. \ref{figure3.1}) such that
 
$ V\left ( D_{1}  \right ) =\left \{ x,~y,~u_{1,2},\cdots,~u_{2,m-1},~u_{2,m},\cdots,~u_{n,m},~u_{1,m}\right \}$ and


$A\left ( D_{1}  \right )=\left \{xu_{1,2},~u_{1,2}y,\cdots,~u_{2,m-1}u_{2,m},\cdots ,~u_{n-1,m}u_{n,m},~u_{n,m}u_{1,m},~\right.\\  \left. u_{1,m}x\right \}$.

$V\left ( D_{2}  \right ) =\left \{ x,~y,~u_{2,1},\cdots,~u_{n-1,2},\cdots,~u_{n-1,m-1},~u_{n-1,m},~u_{n-1,1}\right \}$ and

$ A\left ( D_{2}  \right ) =\left \{ xu_{2,1},~u_{2,1}y,\cdots,~u_{n-2,2}u_{n-1,2},\cdots,~u_{n-1,m-1}u_{n-1,m},~u_{n-1,m}\right.\\  \left.u_{n-1,1},~u_{n-1,1}u_{n,1},~u_{n,1}x\right \} $.

Then we have $2=\min\{\delta^+(D), \delta^-(D)\}\geq
\lambda_2(\overrightarrow{C}_n \Box \overrightarrow{C}_m)\geq 2$.~This completes the proof.
\end{pf}

\begin{pro}\label{p2}
We have $ \lambda_2(\overrightarrow{C}_n\Box \overleftrightarrow{C} _{m} )=3.$
\end{pro}

\begin{figure}[htbp]
\small
\centering
  \includegraphics[width=9cm]{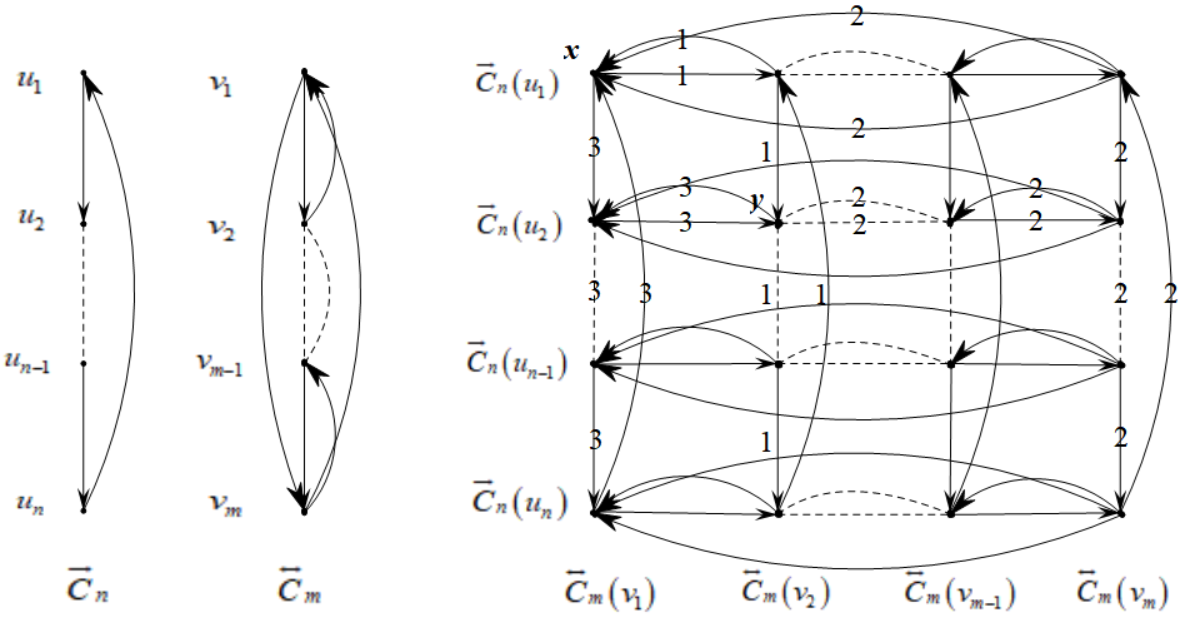}
\caption{Cartesian product of a dicycle and the complete biorientation of a cycle.}
\label{figure3.2}
\end{figure}

\begin{pf}\,   
Let $S=\left \{ x, y \right \} $,~we just consider the case that $ x $,~$ y$ are  neither in the same $ \overrightarrow{C}_n\left ( u_{i}  \right )  $ nor in the same $ \overleftrightarrow{C } _{m} \left ( v_{j}  \right )  $ for some $ 1\leq i \leq n $,~$ 1\leq j \leq m $,~since the arguments for remaining cases are similar.~Without loss of generality,~we may assume that  $x=u_{1,1}$,~$y=u_{2,2}$.~We can  get  three arc-disjoint $S$-strong subgraphs in $\overrightarrow{C}_n\Box \overleftrightarrow{C } _{m} $,~say $ D_{1} $,~$ D_{2} $ and $ D_{3}$ 
(as shown in  Fig. \ref{figure3.2}) such that

$ V\left ( D_{1}  \right ) =\left \{ x,~y,\cdots,~u_{n-1,2},~u_{n,2},~u_{1,2}\right \} $ and


$A\left ( D_{1}  \right )=\left \{ xu_{1,2},~u_{1,2}y,\cdots,~u_{n-1,2}u_{n,2},~u_{n,2}u_{1,2},~u_{1,2}x\right \}$.

$V\left ( D_{2}  \right ) =\left \{ x,~y,~u_{1,m},~u_{2,m},~u_{2,m-1},\cdots,~u_{n-1,m-1},~u_{n,m}\right \}$ and

$ A\left ( D_{2}  \right ) =\left \{ xu_{1,m},~u_{1,m}u_{2,m},~ u_{2,m}u_{2,m-1},\cdots,~u_{2,3}y,~yu_{2,3},\cdots,~u_{2,m-1}\right.\\  \left.u_{2,m},\cdots,~u_{n-1,m}u_{n,m},~u_{n,m}u_{1,m},~u_{1,m}x\right \} $.

$V\left ( D_{3}  \right ) =\left \{ x,~y,~u_{2,1},\cdots,~u_{n-1,1},~u_{n,1}\right \}$ and

$ A\left ( D_{3}  \right ) =\left \{ xu_{2,1},~u_{2,1}y,~ yu_{2,1},\cdots,~u_{n-1,1}u_{n,1},~u_{n,1}x\right \} $.

Then we have
$3=\min\{\delta^+(D),~\delta^-(D)\}~\geq
\lambda_2(\overrightarrow{C}_n \Box \overleftrightarrow{C }_{m} )\geq 3$.~This completes the proof.
\end{pf}

\begin{pro}\label{p4}
We have $ \lambda_2(\overrightarrow{C}_n\Box \overleftrightarrow{T} _{m} )=2.$
\end{pro}
\begin{figure}[htbp]
\small
\centering
  \includegraphics[width=9cm]{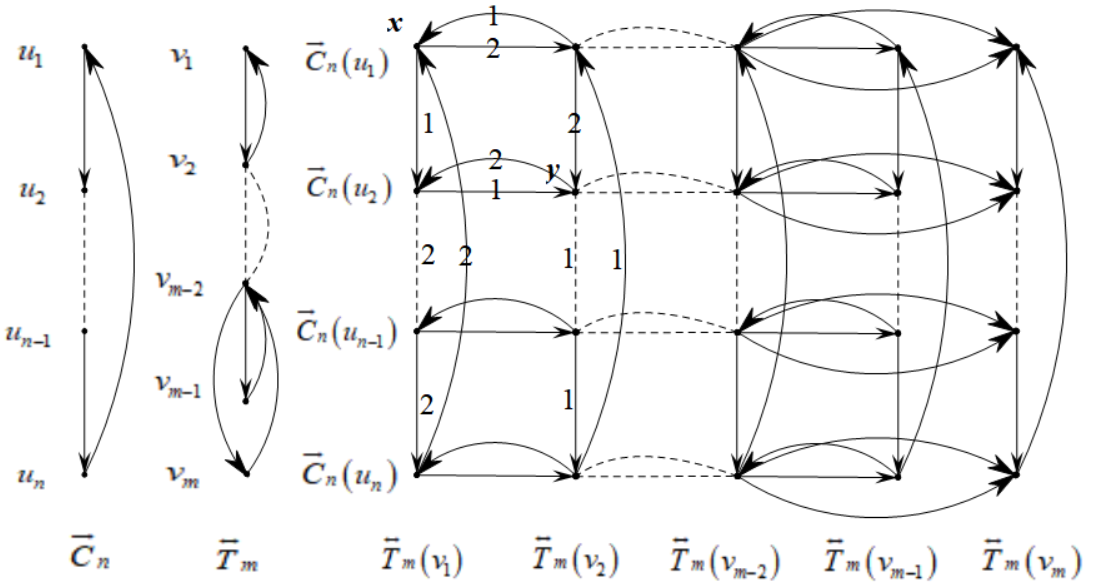}
\caption{Cartesian product of a dicycle and an orientation of a tree.}
\label{figure3.4}
\end{figure}

\begin{pf}\,   
Let $S=\left \{ x, y \right \} $,~we just the case consider that $x$,~$ y$  are neither in the same $ \overrightarrow{C}_n\left ( u_{i}  \right )  $  nor in the same  $ \overleftrightarrow{T } _{m} \left ( v_{j}  \right )  $ for some $ 1\leq i \leq n $,~$ 1\leq j \leq m $,~as the  arguments for the remaining cases are similar.~ Without loss of generality,~we may assume that $ x=u_{1,1}$,~$y=u_{2,2}$.~We can  get two arc-disjoint $S$-strong subgraphs in $\overrightarrow{C}_n\Box \overleftrightarrow{T } _{m}$,~say $ D_{1} $ and $D_{2}$ (as shown in Fig. \ref{figure3.4}) such that

$ V\left ( D_{1}  \right ) =\left \{ x,~y,\cdots,~u_{n-1,2},~u_{n,2},~u_{1,2},~u_{2,1}\right \} $ and


$A\left ( D_{1}  \right )=\left \{ xu_{2,1},~u_{2,1}y,\cdots,~u_{n-1,2}u_{n,2},~u_{n,2}u_{1,2},~u_{1,2}x\right \}$.

$V\left ( D_{2}  \right ) =\left \{x,~y,~u_{1,2},~u_{2,1},\cdots,~u_{n-1,1},~u_{n,1}\right \}$ and

$ A\left ( D_{2}  \right ) =\left \{ xu_{1,2},~u_{1,2}y,~yu_{2,1},\cdots,~u_{n-1,1}u_{n,1},~u_{n,1}x\right \} $.

Then we have $2=\min\{\delta^+(D),~\delta^-(D)\}\geq
\lambda_2(\overrightarrow{C}_n \Box \overleftrightarrow{T } _{m} )\geq 2$.~This completes the proof.
\end{pf}

\begin{pro}\label{p7}
We have $ \lambda_2(\overrightarrow{C}_n\Box \overleftrightarrow{K} _{m} )=m.$
\end{pro}
\begin{figure}[htbp]
\small
\centering
  \includegraphics[width=9cm]{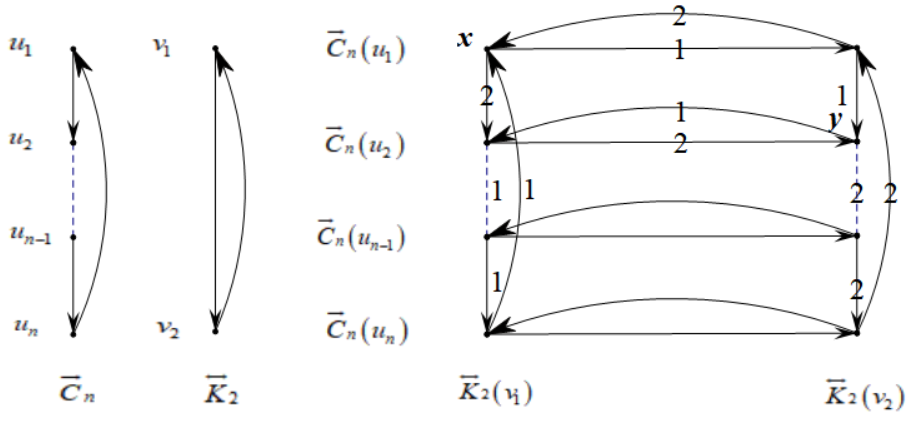}
\caption{Cartesian product of a dicycle and the complete biorientation of $K_2$.}
\label{figure3.7}
\end{figure}

\begin{figure}[htbp]
\small
\centering
  \includegraphics[width=9cm]{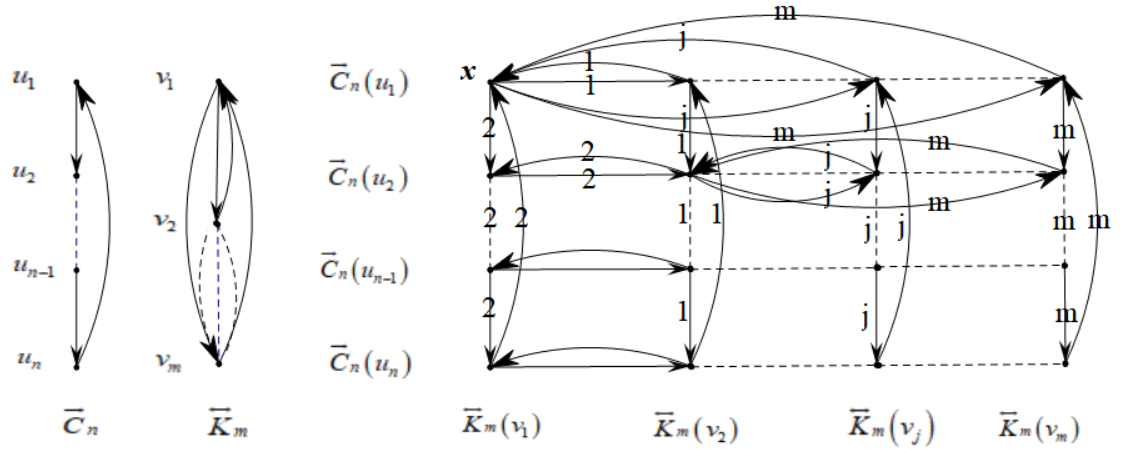}
\caption{Cartesian product of a dicycle and the complete biorientation of $K_m$.}
\label{figure3.7.1}
\end{figure}

\begin{pf}\,   
Let $S=\left \{ x, y \right \} $,~we just consider the case that $ x $,~$ y$  are  neither in the same $ \overrightarrow{C}_n\left ( u_{i}  \right ) $ nor in the same $ \overleftrightarrow{K } _{m} \left ( v_{j}  \right )$ for some $ 1\leq i \leq n $,~$ 1\leq j \leq m $,~as the arguments for the remaining cases are similar.~Without loss of generality,~we may assume that $ x=u_{1,1}$,~$y=u_{2,2}$.

We first show that $ \lambda_2(\overrightarrow{C}_n\Box \overleftrightarrow{K} _{2} )=2 $.~When $m=2$,~we can get two arc-disjoint $S$-strong subgraphs in $\overrightarrow{C}_n\Box \overleftrightarrow{K } _{2} $,~say $ D_{1} $ and $ D_{2} $~(as shown in Fig. \ref{figure3.7}) satisfying:

$ V\left ( D_{1}  \right ) =\left \{x,~y,~u_{1,2},~u_{2,1},\cdots,~u_{n-1,1},~u_{n,1}\right \} $ and


$A\left ( D_{1}  \right )=\left \{ xu_{1,2},~u_{1,2}y,~yu_{2,1},\cdots,~u_{n-1,1}u_{n,1},~u_{n,1}x\right \}$.

$V\left ( D_{2}  \right ) =\left \{x,~y,~u_{1,2},~u_{2,1},\cdots,~u_{n-1,1},~u_{n,1}\right \}$ and

$ A\left ( D_{2}  \right ) =\left \{ xu_{2,1},~u_{2,1}y,\cdots,~u_{n-1,2}u_{n,2},~u_{n,2}u_{1,2},~u_{1,2}x\right \} $.

The propositon is now proved by induction on $m$. Suppose that when $m=k$,~we have $\lambda_2(\overrightarrow{C}_n \Box \overleftrightarrow{K } _{k} )=k$.~We shall show that $\lambda_2(\overrightarrow{C}_n \Box \overleftrightarrow{K } _{k+1} )=k+1 $~when $m=k+1$.Since we can get $k$ arc-disjoint  $S$-strong subgraphs  in $\overrightarrow{C}_n\Box \overleftrightarrow{K } _{k} $ ,~say $ D_{1} $,~$ D_{2} $,$\cdots$,~$ D_{k} $.~When $ m=k+1$,~that is,~the degree of each vertex increases by 2 in $\overleftrightarrow{K } _{k}$,~we can  get  $k+1$ arc-disjoint $S$-strong subgraphs in $\overrightarrow{C}_n\Box \overleftrightarrow{K } _{k+1} $,~say $ D_{1} $,~$ D_{2} $,$\cdots$,~$ D_{k} ,~ D_{k+1} $.~By the symmetry of the complete digraph,~the same conclusion is drawn in the two cases where $x$,~$y$ belong to $\overleftrightarrow{K} _{k+1}$,~and $x$,~$y$ belong to $\overleftrightarrow{K} _{k}$ and $\overleftrightarrow{K} _{k+1}$,~respectively.~From the above argument, the original proposition holds for any positive integer,~we can  get  $m$ arc-disjoint $S$-strong subgraphs in $\overrightarrow{C}_n\Box \overleftrightarrow{K } _{m} $,~say $ D_{1},~D_{2},\cdots,~D_{
j} (2< j\le m),\cdots,~D_{m-1},~D_{m} $
(as shown in Fig. \ref{figure3.7.1}) such that

$ V\left ( D_{1}  \right ) =\left \{x,~y,~u_{1,2},\cdots,~u_{n-1,2},~u_{n,2}\right \} $ and


$A\left ( D_{1}  \right )=\left \{ xu_{1,2},~u_{1,2}y,\cdots,~u_{n-1,2}u_{n,2},~u_{n,2}u_{1,2},~u_{1,2}x\right \}$.

$ V\left ( D_{2}  \right ) =\left \{x,~y,~u_{2,1},\cdots,~u_{n-1,2},~u_{n,2}\right \} $ and


$A\left ( D_{2}  \right )=\left \{ xu_{2,1},~u_{2,1}y,~yu_{2,1},\cdots,~u_{n-1,1}u_{n,1},~u_{n,1}x\right \}$.

$V\left ( D_{j}  \right ) =\left \{x,~y,~u_{1,2},~u_{2,j},\cdots,~u_{n-1,j},~u_{n,j}\right \}$ and

$ A\left (D_{j}  \right ) =\left \{ xu_{1,j},~u_{1,j}u_{2,j},~u_{2,j}y,~yu_{2,j},\cdots,~u_{n-1,j}u_{n,j},~u_{n,j}u_{1,j},~u_{1,j} x\right \} $.

Then we have 
$m=\min\{\delta^+(D), \delta^-(D)\}~\geq
\lambda_2(\overrightarrow{C}_n \Box \overleftrightarrow{K } _{m} )\geq m$.~This completes the proof.
\end{pf}

Since $\lambda_2(\overleftrightarrow{Q})=\lambda(Q)$ for any undirected graph $Q$, using Cartesian product definition, we have 
\begin{equation}\label{eq2}
\lambda_2(\overleftrightarrow{G}\Box \overleftrightarrow{H})= \lambda(G\Box H) 
\end{equation}
for undirected graphs $G$ and $H.$

Propositions \ref{p1}-\ref{p7} and Formulas (\ref{eq2}) and (\ref{eq3}) imply
the following theorem. Indeed, entries in the first row and columns of Table 1 follow from Propositions \ref{p1}-\ref{p7} and all other entries can be easily computed using (\ref{eq2}) and (\ref{eq3}). 

\begin{thm}\label{thmd1}
The following table for the strong
subgraph 2-arc-connectivity of Cartesion products of some digraph classes holds:

\begin{figure}[htbp]
{\tiny
\begin{center}
\renewcommand\arraystretch{3.5}
\begin{tabular}{|p{1.5cm}|p{1.5cm}|p{1.5cm}|p{1.5cm}|p{1.5cm}|}
\hline & $\overrightarrow{C}_m$ & $\overleftrightarrow{C}_m$ &
$\overleftrightarrow{T}_m$ & $\overleftrightarrow{K}_m$
\\\hline

$\overrightarrow{C}_n$ & $2$ & $3$ & $2$ & $m$
\\\hline

$\overleftrightarrow{C}_n$ & $3$ & $4$ & $3$ & $m+1$
\\\hline

$\overleftrightarrow{T}_n$ & $2$ & $3$ & $2$ & $m$
\\\hline

$\overleftrightarrow{K}_n$ & $n$ & $n+1$ & $n$ & $n+m-2$
\\\hline

\end{tabular}
\vspace*{20pt}

\centerline{\normalsize Table $1$. Exact values of $\lambda_2$ for Cartesian products of some digraph classes.}
\end{center}}
\end{figure}
\end{thm}
\vskip 0.5cm

\end{document}